\documentclass{paper}

\usepackage{latexsym,amssymb,amsfonts,amsmath}
\usepackage[pdftex]{graphicx}
\usepackage{natbib}

\title{First moments of the truncated and absolute Student's variates}
\author{Christian P. Robert\\
Universit\'e Paris-Dauphine, IUF, \&~CREST}

\begin{document}

\maketitle

\begin{abstract}
While some of the enclosed already is a well-known derivation, and the remaining may have been obtained
in earlier publications, this note computes the first two moments of
a  Student's variate truncated at zero and of an absolute (or folded) Student's variate.
\end{abstract}

\begin{keywords}Student's t distribution, mean, variance\end{keywords}

\section*{}
This note contains the derivation of both first moments of an absolute Student's variate, also called {\em folded
$t$} in \cite{johnson:kotz:balakrishnan:1994} and \cite{psarakis:panaretos:1990}, who initially derived the
following moments when $\mu=0$, and of
a truncated Student's variate, whose values I could not trace in \cite{johnson:kotz:balakrishnan:1994}.

\section{Absolute Student's variate}
Given a regular Student's $t$ variate
$$
X\sim \mathcal{T}(\nu,\mu,1)\,,
$$
which is the Student's $t$ distribution with location parameter $\mu$ and $\nu>1$ degrees of freedom (and w.l.o.g.~variance $1$), 
we are interested in the first two moments of the transformed rv $|X|$. Since the mean of the 
Student's $t$ distribution only exists for $\nu>1$ we will assume this is the case in the rest of this note.

\begin{figure}
\centerline{\includegraphics[height=.5\textwidth]{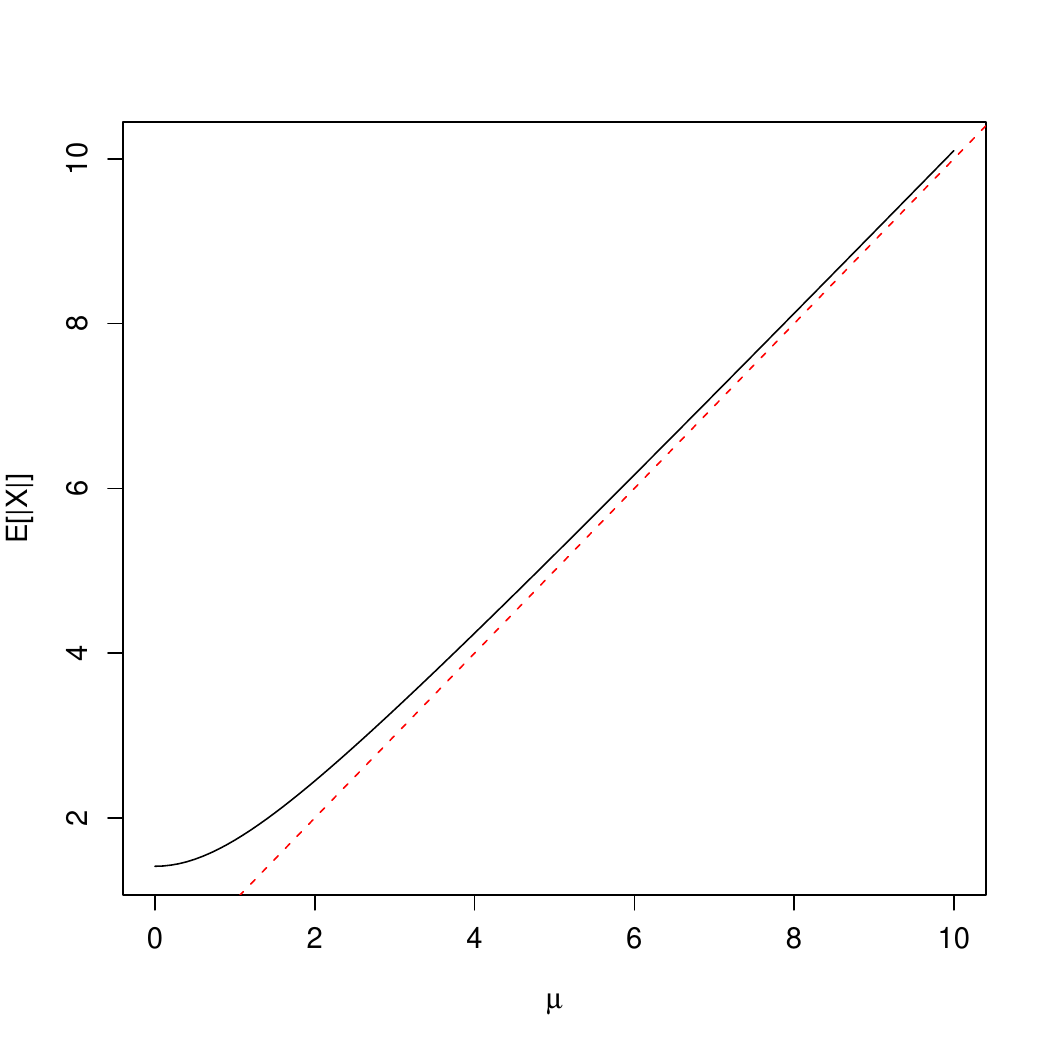}}
\caption{\label{fig:absomean}
Evolution of $\mathbb{E}_{\nu,\mu}[|X|]$ when $\nu=2$ and $\mu$ varies. The dotted line is the diagonal,
asymptote of the expectation.}
\end{figure}

First, we have
$$
\mathbb{E}_{\nu,\mu}[|X|]  =
\int_0^\infty \frac{x}{\left[1+(x-\mu)^2/\nu\right]^{(1+\nu)/2}} C(\nu)\text{d}x -
\int_{-\infty}^0 \frac{x}{\left[1+(x-\mu)^2/\nu\right]^{(1+\nu)/2}} C(\nu)\text{d}x
$$
where $C(\nu)$ is the normalising constant of the usual Student's $t$ density, equal to the 
standard Student's's $t$ density at zero $f(0|\nu)$, i.e.~
$$
C(\nu)=\dfrac{\Gamma((\nu+1)/2)}{\Gamma(\nu/2)\sqrt{\pi\nu}}\,.
$$
Thus,
\begin{align*}
\int_0^\infty \frac{x}{\left[1+(x-\mu)^2/\nu\right]^{(1+\nu)/2}} C(\nu)\text{d}x &= \mu \mathbb{P}_{\nu,\mu}(X>0) \\
- \frac{\nu}{\nu-1}&\int_0^\infty \frac{\text{d}}{\text{d}x}\left\{
	\frac{1}{\left[1+(x-\mu)^2/\nu\right]^{(\nu-1)/2}} 
\right\} C(\nu)\text{d}x\\
&= \mu \mathbb{P}_{\nu,\mu}(X>0) + \frac{C(\nu)\nu/(\nu-1)}{\left[1+\mu^2/\nu\right]^{(\nu-1)/2}}
\end{align*}
and
\begin{align*}
\int_{-\infty}^0 \frac{x}{\left[1+(x-\mu)^2/\nu\right]^{(1+\nu)/2}} C(\nu)\text{d}x &= \mu \mathbb{P}_{\nu,\mu}(X<0) - \\
&\frac{\nu}{\nu-1}\int_{-\infty}^0 \frac{\text{d}}{\text{d}x}
\frac{1}{\left[1+(x-\mu)^2/\nu\right]^{(\nu-1)/2}} C(\nu)\text{d}x\\
&= \mu \mathbb{P}_{\nu,\mu}(X<0) - \frac{C(\nu)\nu/(\nu-1)}{\left[1+\mu^2/\nu\right]^{(\nu-1)/2}}
\end{align*}
so
$$
\mathbb{E}_{\nu,\mu}[|X|]  = 
\mu\left\{2\mathbb{P}_{\nu,\mu}(X>0)-1\right\}+\frac{2 C(\nu)\nu/(\nu-1)}{\left[1+\mu^2/\nu\right]^{(\nu-1)/2}}\,.
$$
In particular, when $\mu=0$,
$$
\mathbb{E}_{\nu,0}[|X|]  = \frac{2 C(\nu)\nu}{\nu-1}\,,
$$
as obtained by \cite{psarakis:panaretos:1990}. Figure \ref{fig:absomean} represents the evolution of
$\mathbb{E}_{\nu,\mu}[|X|]$
when $\nu=2$ and $\mu$ varies. 

For the second moment, $\mathbb{E}_{\nu,\mu}[|X|^2]=\mathbb{E}_{\nu,\mu}[X^2]=\nu/(\nu-2)+\mu^2$ and
therefore $\text{var}_{\nu,\mu}(|X|)$ is easily derived. 
When $\mu=0$, \cite{psarakis:panaretos:1990} (correctly) express this variance as
$$
\text{var}_{\nu,0}(|X|) =
\dfrac{\nu}{\nu-2} - \dfrac{4\nu}{\pi(\nu-1)^2}\,\dfrac{\Gamma((\nu+1)/2)^2}{\Gamma(\nu/2)^2}\,.
$$
(\citealp{johnson:kotz:balakrishnan:1994}, p.403, omitted the minus sign between both terms.)

\section{Truncated Student's variate}
\begin{figure}
\centerline{\includegraphics[height=.5\textwidth]{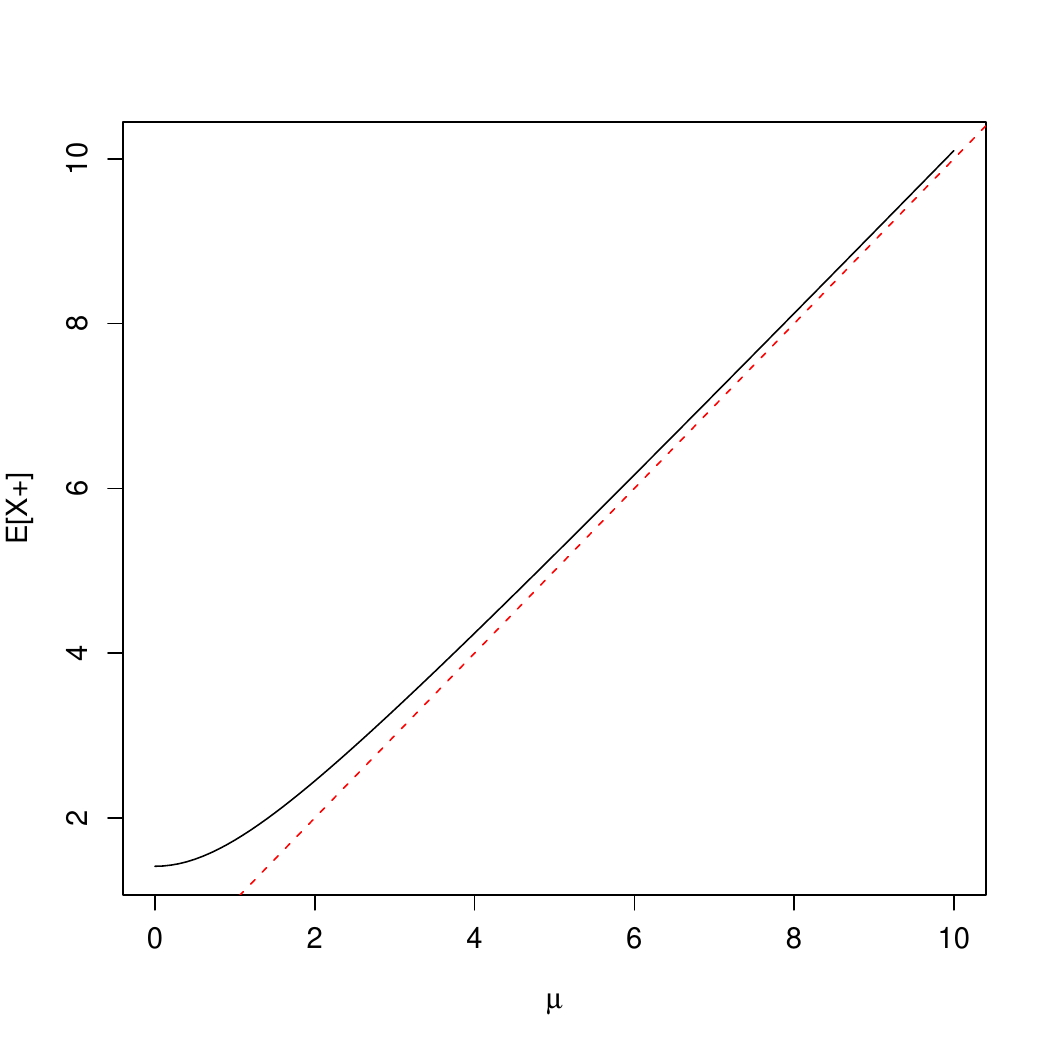}}
\caption{\label{fig:meanplus}
Evolution of $\mathbb{E}_{\nu,\mu}[X_+]$ when $\nu=2$ and $\mu$ varies. The dotted line is the diagonal,
again asymptote of the expectation.}
\end{figure}

We now consider the variate $X_+$ associated with the restriction of the Student's $t$ density to the positive real
line. Its density is provided by
$$
f^+(x|\nu,\mu) = f(x|\nu,\mu) / \mathbb{P}_{\nu,\mu}(X>0)\,.
$$
Therefore, we can derive from the previous calculation of $\mathbb{E}_{\nu,\mu}[|X|]$ that
$$
\mathbb{E}_{\nu,\mu}[X_+] = \mu +
\frac{C(\nu)\nu/(\nu-1)}{\mathbb{P}_{\nu,\mu}(X>0)\,\left[1+\mu^2/\nu\right]^{(\nu-1)/2}}\,.
$$
Figure \ref{fig:meanplus} represents the evolution of $\mathbb{E}_{\nu,\mu}[X_+]$
when $\nu=2$ and $\mu$ varies. 

The second moment is slightly more involved to derive. From ($\nu>2$)
\begin{align*}
\int_0^\infty \dfrac{x^2}{\left[1+(x-\mu)^2/\nu\right]^{(\nu+1)/2}}&\,C(\nu)\,\text{d}x =
\int_0^\infty \dfrac{(x-\mu+\mu)^2}{\left[1+(x-\mu)^2/\nu\right]^{(\nu+1)/2}}\,C(\nu)\,\text{d}x \\
&= \nu \int_0^\infty
\dfrac{\frac{(x-\mu)^2}{\nu}+1-1}{\left[1+(x-\mu)^2/\nu\right]^{(\nu+1)/2}}\,C(\nu)\,\text{d}x \\
&\quad + 2\mu\int_0^\infty \dfrac{x-\mu}{\left[1+(x-\mu)^2/\nu\right]^{(\nu+1)/2}}\,C(\nu)\,\text{d}x\\
&\quad +\mu^2 \mathbb{P}_{\nu,\mu}(X>0)\\
&= (-\nu+\mu^2) \mathbb{P}_{\nu,\mu}(X>0)\\
&\quad + \nu \int_0^\infty \dfrac{C(\nu)}{\left[1+(x-\mu)^2/\nu\right]^{(\nu-1)/2}}\,\text{d}x\\
&\quad + 2\mu \left\{ \mathbb{E}[X_+] - \mu\right\} \mathbb{P}_{\nu,\mu}(X>0)\,,
\end{align*}
we need to compute
\begin{align*}
\int_0^\infty &\dfrac{C(\nu)}{\left[1+(x-\mu)^2/\nu\right]^{(\nu-1)/2}}\,\text{d}x
= \int_0^\infty \dfrac{C(\nu)}{\left[1+\frac{(x-\mu)^2}{(\nu-2)
\frac{\nu}{\nu-2}}\right]^{(\nu-1)/2}}\,\text{d}x \\
&= \dfrac{C(\nu)}{C(\nu-2)}\,\sqrt{\frac{\nu}{\nu-2}}\,\mathbb{P}_{\nu-2,\mu,\sqrt{\nu/(\nu-2)}}(X>0) \,,
\end{align*}
where the third index denotes the scale factor of a Student't distribution,
i.e.~the positivity probability is computed for a $\mathcal{T}(\nu-2,\mu,\sqrt{\nu/(\nu-2)})$
distribution.  Therefore,
$$
\mathbb{E}_{\nu,\mu}[X_+^2] = -(\nu+\mu^2) + \dfrac{\nu
C(\nu)}{C(\nu-2)}\,\sqrt{\frac{\nu}{\nu-2}}\,
\dfrac{\mathbb{P}_{\nu-2,\mu,\sqrt{\nu/(\nu-2)}}(X>0)}{\mathbb{P}_{\nu,\mu}(X>0)}+2\mu \mathbb{E}_{\nu,\mu}[X_+]\,.
$$
The corresponding $\text{var}_{\nu,\mu}(X_+)$ is straightforward to derive.

\end{document}